\documentclass[11pt]{amsart}
\usepackage{amssymb, latexsym}

\numberwithin{equation}{section}
\newtheorem{Thm}[equation]{Theorem}
\newtheorem{Prop}[equation]{Proposition}
\newtheorem{Cor}[equation]{Corollary}
\newtheorem{Lem}[equation]{Lemma}

\theoremstyle{definition}
\newtheorem{Def}[equation]{Definition}

\newtheorem{Rmk}[equation]{Remark}

\begin{document}

\title [Iwahori-Hecke Algebras of ${SL}_2$]
{Iwahori-Hecke Algebras of ${SL}_2$\\
over $2$-dimensional Local Fields}
\author[Kyu-Hwan Lee]{Kyu-Hwan Lee$^{\star}$}
\thanks{$^{\star}$ supported in part by EPSRC grant on zeta
functions and in part by KOSEF Grant \# 
R01-2003-000-10012-0.}
\address{Department of Mathematics,
University of Connecticut, Storrs, CT 06269-3009, U.S.A.}
\email{khlee@math.uconn.edu}

\begin{abstract}
In this paper we construct an analogue of Iwahori-Hecke algebras
of ${SL}_2$ over $2$-dimensional local fields. After considering coset decompositions of double cosets of a Iwahori subgroup, we define a
convolution product on the space of certain functions on $SL_2$, and prove that the product is 
well-defined, obtaining a Hecke algebra. Then we investigate the structure of the Hecke algebra. We determine the center of the Hecke algebra and consider Iwahori-Matsumoto type relations.
\end{abstract}

\maketitle

\section*{Introduction}

Hecke algebras play important roles
in the representation theory of $p$-adic groups.
There are two important classes of Hecke algebras.
One is the spherical Hecke algebra attached
to a maximal compact open subgroup, and the other is
the Iwahori-Hecke algebra attached to a Iwahori
subgroup. The spherical Hecke algebra is isomorphic
to the center of the corresponding Iwahori-Hecke algebra.
In the theory of higher dimensional local fields \cite{LF}, a $p$-adic field is a $1$-dimensional local field. So the theory of $p$-adic groups and their Hecke algebras is over $1$-dimensional local fields.

Recently, the representation theory of algebraic groups over
$2$-dimensional local fields was initiated by the works of
Kapranov, Kazhdan and Gaitsgory \cite{Ka,GK1,GK2,GK3}. 
In their
development, Cherednik's double affine Hecke algebras (\cite{Ch}) appear as an
analogue of Iwahori-Hecke algebras.
The common feature of the works mentioned above is the use of rank one integral structure of a $2$-dimensional local field. But there is also a rank two integral structure in a 
$2$-dimensional local field, and it is important in the arithmetic theory to use the rank two integral structure \cite{LF}; we refer the reader to Fesenko's article \cite{Fe1} and to the references there for recent developments.

In their paper \cite{KL}, Kim and Lee constructed an analogue
of spherical Hecke algebras of ${SL}_2$ over $2$-
dimensional local fields using rank two integral structure. They also established a Satake isomorphism using
Fesenko's $\mathbb{R}((X))$-valued measure defined in \cite{Fe} (See also \cite{Fe1}).
In particular, the algebra is proved to be commutative. A similar
result is expected in the case of ${GL}_n$ \cite{KL1}  and, eventually,
in the case of reductive algebraic groups.

In this paper, we construct an analogue of Iwahori-Hecke algebras
of ${SL}_2$ over $2$-dimensional local fields coming from rank two integral structure. Our basic approach will be similar to that of \cite{KL}. More precisely, an element of
the algebra is an infinite linear combination of characteristic
functions of double cosets satisfying certain conditions. We
define a convolution product of two characteristic functions,
using coset decompositions of double cosets of an Iwahori subgroup and imposing a restriction on the support of the
product. Since a double coset of the Iwahori subgroup is an
uncountable union of cosets in general, it is necessary to prove
that the convolution product is well-defined. This will be done in
the first part of the paper.

In the second part, we study the structure of the Hecke algebra. It has a natural $\mathbb{Z}$-grading and contains the affine Hecke algebra of $SL_2$ as a subalgebra. We will find a big commutative subalgebra, and determine its structure completely. Surprisingly, it is different from the group algebra of double cocharacters. We will also calculate the center of the Hecke algebra, and it turns out to be the same as the center of the affine Hecke algebra of $SL_2$.
The classical Iwahori-Hecke algebra has a well-known
presentation due to Iwahori and Matsumoto \cite{I,IM}.
The relations can be understood as deformations of
Coxeter relations of the affine Weyl group. The corresponding
Weyl group of a reductive algebraic
group over a $2$-dimensional
local field is not a Coxeter group.
In ${SL}_2$ case, A.N. Parshin obtained an explicit
presentation of the (double affine) Weyl group \cite{Pa}.
We will investigate Iwahori-Matsumoto
type relations of the Hecke algebra in light of Parshin's presentation.

Now that we have spherical Hecke algebras and Iwahori-Hecke
algebras of ${SL}_2$ over $2$-dimensional local fields with respect to rank two integral structure, next natural steps would be considering representations of $SL_2$ over $2$-dimensional local fields, constructing the Hecke
algebras for more general reductive algebraic groups, and understanding these algebras in connection with the works  of Kapranov, Kazhdan and Gaitsgory mentioned earlier. It seems that another interesting approach to the representation theory over two-dimensional local fields could be obtained from the work of Hrushovski and Kazhdan \cite{HK}, in which they developed a theory of motivic integration. Actually, in the appendix of the paper \cite{HK}, given by Avni, an Iwahori-Hecke algebra of $SL_2$ is constructed using motivic integration. It would be nice if one can find any connection of it to the constructions of this paper.

There are three sections and an appendix in this paper. In the first
section, we fix notations and recall the Bruhat decomposition. The
next section is devoted to the construction of the Hecke algebra.
We will show that the convolution product is well-defined. In the
third section, we study the structure of the Hecke algebra, determining the center of the algebra and finding Iwahori-Matsumoto type relations. In the appendix, we present a complete set of
formulas for convolution products of characteristic functions. These formulas will be essentially used for many calculations in this paper. 

\subsection*{Acknowledgments} \quad
I would like to thank
Henry Kim and Ivan Fesenko for their encouragement and
useful comments during preparation of this paper.

\vskip 1 cm

\section{Bruhat Decomposition}

In this section, we fix notations and collect some results on double cosets and coset decompositions we will use later.
We assume that the reader is familiar with basic definitions in the theory of 
$2$-dimensional local fields, which can be found in \cite{Zu}. 

\vskip 5mm

Let $F(=F_2)$ be a two dimensional local field with the first
residue field $F_1$ and the last residue field
$F_0(=\mathbb{F}_q)$ of $q$ elements. We fix a discrete valuation
$v: F^{\times}
\rightarrow \mathbb{Z}^2$ of rank two.  Recall that $\mathbb{Z}^2$
is endowed with the lexicographic ordering from the right.  Let
$t_1$ and $t_2$ be local parameters with respect to the valuation $v$.
We have the ring $O$ 
of integers of $F$ with respect to the rank-two valuation $v$. There is a natural projection
$p : O \rightarrow O/t_1 O = F_0$. Note that the ring $O$ is different from the ring $O_{21}$ of integers of $F$ with respect to the rank-one valuation $v_{21}: F^{\times}
\rightarrow \mathbb{Z}$.

\vskip 5mm

Let ${G}$ be a connected split semisimple algebraic group
defined over $\mathbb{Z}$. We fix a maximal torus ${T}$ and
a Borel subgroup ${B}$ such that ${T} \subset
{B} \subset {G}$, and we have $W_0=N_{G}({T})/{{T}}$, the Weyl group of ${G}$. We write $G=G(F)$ and consider $I = \{ x \in {G}(O) :
p(x) \in
{B}(F_0) \}$, the {\em double Iwahori subgroup}, and $W = N_G(T)/T(O)$, the {\em double affine Weyl group}, and obtain the following decomposition.

\begin{Prop} \cite{Ka,Pa}\label{Bruhat} We have
   \[ G = \coprod_{w
\in W} I w I , \] and the resulting identification $I \backslash
G/I \rightarrow W$ is independent of the choice of representatives
of elements of $W$.
\end{Prop}

From now on, in the rest of this paper,
we assume that ${G}={SL}_2$ and $B$ is the subgroup of upper triangular matirices.
The following lemma gives explicit formulas for the
decomposition in Proposition \ref{Bruhat}. The proof is straightforward, so we omit it.

\begin{Lem} \label{lem-explicit} Assume that
$x=\begin{pmatrix} a & b \\ c & d \end{pmatrix} \in G$.
\begin{enumerate} 
\item If $v(a) \le v(b)$ and $v(a)<v(c)$, then $ x \in I
\begin{pmatrix} a & 0 \\ 0 & a^{-1} \end{pmatrix} I$.
\vskip 5 mm
\item If $v(b)<v(a)$ and $v(b) < v(d)$, then $x \in I
\begin{pmatrix} 0 & b \\ -b^{-1} & 0 \end{pmatrix} I$.
\vskip 5 mm
\item If $v(c) \le v(a)$ and $v(c) \le v(d)$, then $x \in I
\begin{pmatrix} 0 & -c^{-1} \\ c & 0 \end{pmatrix} I$.
\vskip 5 mm
\item If $v(d) \le v(b)$ and $v(d)<v(c)$, then $x \in I
\begin{pmatrix} d^{-1} & 0 \\ 0 & d \end{pmatrix} I$.
\end{enumerate}
\end{Lem}

\vskip 5mm

We denote by $C^{(1)}_{i,j}$ and
$C^{(2)}_{i,j}$, $(i,j) \in \mathbb{Z}^2$,
the double cosets
\[ I \begin{pmatrix}
t_1^i t_2^j & 0 \\ 0 & t_1^{-i} t_2^{-j} \end{pmatrix} I
\quad \text{and} \quad I \begin{pmatrix} 0 &
t_1^i t_2^j  \\ - t_1^{-i} t_2^{-j} & 0
\end{pmatrix} I, \quad \text{respectively}.  \]
In the following lemma, we obtain complete sets of coset
representatives in the decomposition of double cosets of the subgroup
$I$ into right cosets.

\begin{Lem} \label{lem-rep}
We have $C^{(a)}_{i,j}= \coprod_z Iz$ where the disjoint
union is over $z$ in the following list.
\begin{enumerate}
\item If $a=1$ and $(i,j) \ge (0,0)$, then \[
z=\begin{pmatrix}
t_1^i t_2^j & 0 \\ 0 & t_1^{-i} t_2^{-j} \end{pmatrix},
\quad \begin{pmatrix}
t_1^i t_2^j & 0 \\ t_1^{-i+k} t_2^{-j+l}u
 & t_1^{-i} t_2^{-j} \end{pmatrix} \text{ for }
(1,0) \le (k,l) < (2i+1,2j), \] where $u \in O^{\times}$
 are units belonging to a fixed set of representatives
 of $O/t_1^{2i-k+1}t_2^{2j-l}O$.

\item If $a=1$ and $(i,j) < (0,0)$, then \[
z=\begin{pmatrix}
t_1^i t_2^j & 0 \\ 0 & t_1^{-i} t_2^{-j} \end{pmatrix},
\quad \begin{pmatrix}
t_1^i t_2^j & t_1^{i+k} t_2^{j+l}u \\ 0
 & t_1^{-i} t_2^{-j} \end{pmatrix} \text{ for }
(0,0) \le (k,l) < (-2i,-2j), \] where $u \in O^{\times}$
 are units belonging to a fixed set of representatives
 of $O/t_1^{-2i-k}t_2^{-2j-l}O$.

\item If $a=2$ and $(i,j) \ge (0,0)$, then \[
z=\begin{pmatrix}
0 & t_1^i t_2^j \\ -t_1^{-i} t_2^{-j} & 0 \end{pmatrix},
\quad \begin{pmatrix}
0 & t_1^i t_2^j \\ -t_1^{-i} t_2^{-j}
 & - t_1^{-i+k} t_2^{-j+l}u \end{pmatrix} \text{ for }
(0,0) \le (k,l) < (2i+1,2j), \] where $u \in O^{\times}$
 are units belonging to a fixed set of representatives
 of $O/t_1^{2i-k+1}t_2^{2j-l}O$.

\item If $a=2$ and $(i,j) < (0,0)$, then \[
z=\begin{pmatrix}
0 & t_1^i t_2^j \\ -t_1^{-i} t_2^{-j} & 0 \end{pmatrix},
\quad \begin{pmatrix}
t_1^{i+k} t_2^{j+l}u & t_1^i t_2^j \\ -t_1^{-i} t_2^{-j}
 & 0 \end{pmatrix} \text{ for }
(1,0) \le (k,l) < (-2i,-2j), \] where $u \in O^{\times}$
 are units belonging to a fixed set of representatives
 of $O/t_1^{-2i-k}t_2^{-2j-l}O$.
\end{enumerate}
\end{Lem}

\begin{proof}
Since the other cases are similar, we only prove
the part (1). Consider
\[ z = \begin{pmatrix}
t_1^i t_2^j & 0 \\ 0 & t_1^{-i} t_2^{-j} \end{pmatrix}
\begin{pmatrix}
a & b \\ c & d \end{pmatrix}, \qquad
z'=\begin{pmatrix}
t_1^i t_2^j & 0 \\ 0 & t_1^{-i} t_2^{-j} \end{pmatrix}
\begin{pmatrix}
a' & b' \\ c' & d'  \end{pmatrix}\] where
$\begin{pmatrix}
a & b \\ c & d \end{pmatrix},
  \begin{pmatrix}
a' & b' \\ c' & d'  \end{pmatrix} \in I$. We see that
the condition $Iz=Iz'$ is equivalent to \[c'd -cd'
\in t_1^{2i+1} t_2^{2j}O .\] We write $(c,d) \sim
(c',d')$ if $  c'd -cd'
\in t_1^{2i+1} t_2^{2j}O$.
Note that if $\begin{pmatrix} a & b \\
c & d \end{pmatrix} \in I$ then $c \in t_1 O$ and
$d$ is a unit. Let
$C$ be the set of pairs $(c,d) \in O^2$ such that
$c \in t_1 O$ and
$d$ is a unit. Then $\sim$ is an
equivalence relation on $C$. In order to determine different
cosets, we need only to determine a set of representatives of the
equivalence relation $\sim$, which turn out
to be \[ (0,1)\quad \text{and} \quad
 (t_1^k t_2^lu, 1)  \text{ for }
(1,0) \le (k,l) < (2i+1,2j), \] where $u \in O^{\times}$
 are units belonging to a fixed set of representatives
 of $O/t_1^{2i-k+1}t_2^{2j-l}O$. These
yield the elements $z$ in the part (1).
\end{proof}

\begin{Rmk} \label{rmk-inf}
The disjoint union $C^{(a)}_{i,j}= \coprod_z Iz$ is an uncountable union unless $j=0$. The same
is true for the double cosets
of $K={SL}_2(O)$; see \cite{KL}.
\end{Rmk}

\vskip 1 cm

\section{Iwahori-Hecke algebras}

In this section, we define the convolution product of two
characteristic functions of double cosets of the subgroup $I$, and we show that the
product is well-defined. Then we construct an analogue of the 
Iwahori-Hecke algebra of ${SL}_2$.

\vskip 5mm

We fix a set of representatives $R$ of the double affine Weyl group $W$ to be \begin{equation} \label{representatives} R=
\left \{ \eta^{(1)}_{i,j} :=
\begin{pmatrix}
t_1^i t_2^j & 0 \\ 0 & t_1^{-i} t_2^{-j} \end{pmatrix},
\quad   \eta^{(2)}_{i,j} :=\begin{pmatrix} 0 &
t_1^i t_2^j  \\ - t_1^{-i} t_2^{-j} & 0
\end{pmatrix} : \quad (i,j) \in \mathbb{Z}^2 \right \}. \end{equation}
We define a map $\eta: G \rightarrow R$ so that $x \in I\eta(x)I$ for each $x \in G$. That is, we assign to an element $x$ of $G$ its representative $\eta(x) \in R$ in the decomposition $G = \coprod_{w
\in W} I w I$ of Proposition \ref{Bruhat}.

We put
\begin{equation} \label{eqn-cell} C_j = \coprod_{m \in \mathbb{Z}} C^{(1)}_{m,j}
\bigcup C^{(2)}_{m,j}, \qquad  j \in \mathbb{Z}.\end{equation}

We denote by $\chi^{(1)}_{i,j}$ and $\chi^{(2)}_{i,j}$
the characteristic
functions of the double cosets $C^{(1)}_{i,j}$ and
$C^{(2)}_{i,j}$, respectively.  We will consider the following types of functions
\begin{equation} \label{eqn-type}
  \sum_{r \le i} c_r
\chi^{(a)}_{r,j} \ (j>0), \quad  \sum_{r \ge i} c_r
\chi^{(a)}_{r,j} \ ( j<0), \quad  \text{and} \quad \sum_{i \le r \le i'} c_r
\chi^{(a)}_{r,0}  \end{equation}  for $i, i', j
 \in \mathbb{Z}$, $a=1,2$ and $c_r \in \mathbb{C}$.
Now we define the convolution product
of two characteristic functions.

\begin{Def}   For $a,b =1,2$ and $(i,j), (k,l)
\in \mathbb{Z}^2$, we define
\begin{equation} \label{convolution-I} \left ( \chi^{(a)}_{i,j}*\chi^{(b)}_{k,l} \right )
(x) = \left \{ \begin{array}{ll} q^{-1}\sum_z
\chi^{(a)}_{i,j}
\left ( \eta(x) z^{-1} \right ) & \text{if } x \in C_{j+l}
\text{ and } jl \ge 0,\\ 0
& \text{otherwise}, \end{array} \right .
 \end{equation} where the sum is over
the representatives $z$ of the decomposition $
C^{(b)}_{k,l} = \coprod_z Iz $.
\end{Def}

\begin{Rmk}
One can construct a certain invariant $\Bbb R((X))$-valued measure
$d\gamma$ on $G$. Then we could define the convolution
product of two functions $f$ and $g$ on $G$ by
$$
(f*g)(x) = \int_G f(xy^{-1}) g(y) d\gamma (y) . 
$$
The definition of the convolution product given above
is derived from this formula.
\end{Rmk}

Since the cardinality of the set of the representatives $z$
in the union is uncountable
(Remark  \ref{rmk-inf}) in general, we need to prove the following.

\begin{Thm} \label{thm-product}
The convolution product 
$\chi^{(a)}_{i,j}*\chi^{(b)}_{k,l}$ yields a well-defined
function  of one of the types in (\ref{eqn-type}) for any
$a,b =1,2$ and $(i,j), (k,l)
\in \mathbb{Z}^2$.
\end{Thm}

\begin{proof}

We write
\[  \chi^{(a)}_{i,j}*\chi^{(b)}_{k,l} = \sum_m
c_m^{(1)} \chi^{(1)}_{m, j+l} + \sum_m c_m^{(2)} \chi^{(2)}_{m,
j+l}, \qquad c_m^{(1)}, c_m^{(2)} \in \mathbb{C} , \ m \in
\mathbb{Z}. \]

\begin{enumerate}
\item Assume that $b=1$, $j \ge 0$ and $(k,l) \ge (0,0)$. By Lemma
\ref{lem-rep}, we need only to consider $z=
\begin{pmatrix}
t_1^k t_2^l & 0 \\ t_1^{-k+k'} t_2^{-l+l'}u
 & t_1^{-k} t_2^{-l} \end{pmatrix} \text{ for }
(1,0) \le (k',l') < (2k+1,2l)$, where $u \in O^{\times}$
 are units belonging to a fixed set of representatives
 of $O/t_1^{2k-k'+1}t_2^{2l-l'}O$.

\begin{enumerate}
\item If $\eta(x)=\eta^{(1)}_{m,j+l}$, then \[ \eta(x) z^{-1}
= \begin{pmatrix} t_1^{m-k}t_2^j & 0 \\
- t_1^{-m-k+k'}t_2^{-j-2l+l'}u & t_1^{-m+k} t_2^{-j}
\end{pmatrix} .\]

\begin{enumerate}
\item If $a=1$, then we have $m-k=i$, and either
$(-m-k+k', -j-2l+l')
>(m-k, j)$ or
$(-m-k+k', -j-2l+l')
>(-m+k, -j)$ by Lemma \ref{lem-explicit}. The first
case gives $(k',l')>(2m, 2j+2l)=(2i+2k, 2j+2l)$, and
so $j=0$, $l'=2l$ and $2i+2k<k'<2k+1$. Counting
the number of possible representative units in
$O/t_1^{2k-k'+1}O$,
we see that $c_{i+k}^{(1)}$
is finite and $c_m^{(1)}=0$ for all $m \neq i+k$.
The second case gives $(k',l')>(2k,2l)$, but
$(k',l')<(2k+1,2l)$ from the assumption, and so
$c_m^{(1)}=0$ for all $m \in \mathbb Z$.

\item If $a=2$, then we must have $l'=2l$ and
$-m-k+k'=-i$. Then $k' < 2k+1$ implies $m < i+k+1$,
and counting
the number of possible representative units in
$O/t_1^{2k-k'+1}O$, we see that
$c_m^{(2)}$ is finite for each $m$.

\end{enumerate}

\item If $\eta(x)=\eta^{(2)}_{m,j+l}$, then  \[ \eta(x) z^{-1}
= \begin{pmatrix} -t_1^{m-k+k'}t_2^{j+l'} u &
t_1^{m+k} t_2^{j+2l} \\
- t_1^{-m-k}t_2^{-j-2l} & 0
\end{pmatrix} .\]

\begin{enumerate}
\item If $a=1$, then $m-k+k'=i$ and $l'=0$, and
it follows from
Lemma \ref{lem-explicit} that
$j \le -l$. Since $j \ge 0$ and $l \ge 0$, we have $j=l=l'=0$.
The condition $1 \le k' < 2k+1$ gives
$i-k-1 < m \le i+k-1$ and
$c_m^{(2)}$ is finite for each $m$.

\item If $a=2$, then $l=l'=0$, $m+k=i$ and
 $1 \le k' < 2k+1$.
Thus $c_{i-k}^{(2)}$ is finite and $c_m^{(2)}=0$ for
$m \neq i-k$.

\end{enumerate}

\end{enumerate}

\item Assume that $b=1$, $j \le 0$ and $(k,l) <(0,0)$. By Lemma
\ref{lem-rep}, we need only to consider $z=
\begin{pmatrix}
t_1^k t_2^l & t_1^{k+k'} t_2^{l+l'}u \\ 0
 & t_1^{-k} t_2^{-l} \end{pmatrix} \text{ for }
(0,0) \le (k',l') < (-2k,-2l)$, where $u \in O^{\times}$
 are units belonging to a fixed set of representatives
 of $O/t_1^{-2k-k'}t_2^{-2l-l'}O$.

\begin{enumerate}
\item If $\eta(x)=\eta^{(1)}_{m,j+l}$, then \[ \eta(x) z^{-1}
= \begin{pmatrix} t_1^{m-k}t_2^j & - t_1^{m+k+k'}
t_2^{j+2l+l'}u \\
0 & t_1^{-m+k} t_2^{-j}
\end{pmatrix} .\]

\begin{enumerate}
\item If $a=1$, then we have $m-k=i$, and either
$(m-k,j) \le
(m+k+k', j+2l+l')$  or $(-m+k,-j) \le
(m+k+k', j+2l+l')$ by Lemma \ref{lem-explicit}.
The first case yields $(k',l') \ge (-2k,-2l)$,
but $(k',l')<(-2k,-2l)$ from the assumption.
Thus $c_m^{(1)}=0$ for all $m \in \mathbb Z$. The
second case gives $(k',l') \ge (-2m, -2j-2l)=
(-2i-2k, -2j-2l)$, and we must have $j=0$, $l'=-2l$
and $-2i -2k \le k' < -2k$. Thus $c^{(1)}_{i+k}$
is finite and $c^{(1)}_m=0$ for all $m \neq i+k$.

\item If $a=2$, then  we have $l'=-2l$ and $m+k+k'=i$.
The condition $k'<-2k$
leads to $i+k<m$,
and $c_m^{(2)}$ is finite for each $m$.

\end{enumerate}

\item If $\eta(x)=\eta^{(2)}_{m,j+l}$, then  \[ \eta(x) z^{-1}
= \begin{pmatrix} 0 & t_1^{m+k}t_2^{j+2l} \\
- t_1^{-m-k}t_2^{-j-2l} & t_1^{-m+k+k'} t_2^{-j+l'}u
\end{pmatrix} .\]

\begin{enumerate}
\item If $a=1$, then $-m+k+k'=-i$ and $l'=0$, and
we have $-j \le l$ by Lemma \ref{lem-explicit}. Since $j \ge 0$ and $l \ge 0$, we have  $j=l=l'=0$.
The condition
$0 \le k' < -2k$ yields $i+k \le m < i-k$ and
$c_m^{(2)}$ is finite for each $m$.

\item If $a=2$, then $l=l'=0$, $m+k=i$ and
$0 \le k' < -2k$. Thus $c_{i-k}^{(2)}$ is finite
and $c_m^{(2)}=0$ for
$m \neq i-k$.

\end{enumerate}

\end{enumerate}

\item Assume that $b=2$, $j \ge 0$ and $(k,l) \ge (0,0)$. By Lemma
\ref{lem-rep}, we need only to consider $z=
\begin{pmatrix}
0 & t_1^k t_2^l \\ -t_1^{-k} t_2^{-l}
 & - t_1^{-k+k'} t_2^{-l+l'}u \end{pmatrix} \text{ for }
(0,0) \le (k',l') < (2k+1,2l)$, where $u \in O^{\times}$
 are units belonging to a fixed set of representatives
 of $O/t_1^{2k-k'+1}t_2^{2l-l'}O$.

\begin{enumerate}
\item If $\eta(x)=\eta^{(1)}_{m,j+l}$, then  \[ \eta(x) z^{-1}
= \begin{pmatrix} -t_1^{m-k+k'}t_2^{j+l'}u &
- t_1^{m+k} t_2^{j+2l} \\
t_1^{-m-k}t_2^{-j-2l} & 0
\end{pmatrix} .\]

\begin{enumerate}
\item If $a=1$, then it is similar to (1)(b)(i).

\item If $a=2$, then it is similar to (1)(b)(ii).

\end{enumerate}

\item If $\eta(x)=\eta^{(2)}_{m,j+l}$, then  \[ \eta(x) z^{-1}
= \begin{pmatrix} t_1^{m-k}t_2^j & 0 \\
t_1^{-m-k+k'}t_2^{-j-2l+l'}u & t_1^{-m+k} t_2^{-j}
\end{pmatrix} .\]

\begin{enumerate}
\item If $a=1$, then  it is similar to (1)(a)(i).

\item If $a=2$, then it is similar to (1)(a)(ii).

\end{enumerate}

\end{enumerate}

\item Assume that $b=2$, $j \le 0$ and $(k,l) <(0,0)$. By Lemma
\ref{lem-rep}, we need only to consider $z=
\begin{pmatrix}
t_1^{k+k'} t_2^{l+l'}u & t_1^k t_2^l \\ -t_1^{-k} t_2^{-l}
 & 0 \end{pmatrix} \text{ for }
(1,0) \le (k',l') < (-2k,-2l)$, where $u \in O^{\times}$
 are units belonging to a fixed set of representatives
 of $O/t_1^{-2k-k'}t_2^{-2l-l'}O$.

\begin{enumerate}
\item If $\eta(x)=\eta^{(1)}_{m,j+l}$, then   \[ \eta(x) z^{-1}
= \begin{pmatrix}  0 & -t_1^{m+k}t_2^{j+2l} \\
t_1^{-m-k}t_2^{-j-2l} & t_1^{-m+k+k'} t_2^{-j+l'}u
\end{pmatrix} .\]

\begin{enumerate}
\item If $a=1$, then  it is similar to (2)(b)(i).

\item If $a=2$, then  it is similar to (2)(b)(ii).

\end{enumerate}

\item If $\eta(x)=\eta^{(2)}_{m,j+l}$, then  \[ \eta(x) z^{-1}
= \begin{pmatrix} t_1^{m-k}t_2^j & t_1^{m+k+k'}
t_2^{j+2l+l'}u \\
0 &  t_1^{-m+k} t_2^{-j}
\end{pmatrix} .\]

\begin{enumerate}
\item If $a=1$, then it is similar to (2)(a)(i).

\item If $a=2$, then it is similar to (2)(a)(ii).

\end{enumerate}

\end{enumerate}

\end{enumerate}

\end{proof}

A complete set of formulas for the
convolution product
$\chi^{(a)}_{i,j}*\chi^{(b)}_{k,l}$ can be found
in Appendix, and we obtain the following result.

\begin{Cor} \label{multiplication}
\begin{equation} \label{eqn-formula} \chi^{(a)}_{i,j}*\chi^{(b)}_{k,l}= \begin{cases} \displaystyle{\sum_{m \le i+k}} c_m \chi^{(b)}_{m, j+l}, \quad c_m \in \mathbb{C}, & \text{ if } a=2, \ j>0 \text{ and } l>0, \\ \displaystyle{\sum_{m > i+k}} c_m \chi^{(b)}_{m, j+l}, \quad  c_m \in \mathbb{C}, & \text{ if } a=2, \ j<0 \text{ and } l<0, \\ \text{\rm a finite sum} &\text{\rm \ otherwise.}\end{cases} \end{equation} 
\end{Cor}

We denote by $\mathcal{H}(G,I)$ the $\mathbb{C}$-vector
space generated by the functions of types in
(\ref{eqn-type}). 
We linearly extend the convolution product $*$ defined in (\ref{convolution-I}) to the whole space $\mathcal{H}(G,I)$. It follows from Corollary \ref{multiplication} that it is well-defined. Thus we have obtained a $\mathbb{C}$-algebra structure on the
space $\mathcal{H}(G,I)$.

\begin{Def}
The $\mathbb{C}$-algebra $\mathcal{H}(G,I)$ will be called the {\em Iwahori-Hecke algebra} of $G(=SL_2)$.
\end{Def}

It can be easily checked that $\iota := q \chi^{(1)}_{0,0}$ is the identity element of the algebra $\mathcal{H}(G,I)$. Furthermore, we have:

\begin{Prop}
The algebra $\mathcal{H}(G,I)$ is an associative algebra.
\end{Prop}

\begin{proof}
Assume that $\alpha= \chi^{(a)}_{i,j}$, $\beta= \chi^{(b)}_{k,l}$ and  $\gamma= \chi^{(c)}_{m,n}$. We fix the sets of representatives $\{ z_1 \}$, $\{ z_2\}$ and  $\{ z_3 \}$ in the decompositions $C^{(a)}_{i,j}= \coprod_{z_1} I z_1$,  $C^{(b)}_{k,l}= \coprod_{z_2} I z_2$ and $C^{(c)}_{m,n}= \coprod_{z_3} I z_3$, respectively. We write \[ \alpha * \beta = \sum_\sigma c(\alpha, \beta; \sigma) \sigma, \text{ where } \sigma = \chi^{(d)}_{p, j+l}, \ d=1, 2, \ p \in \mathbb{Z}. \] Then \[ c(\alpha, \beta; \sigma) = \mathrm{Card} \{ z_2 : \eta^{(d)}_{p,j+1} z_2^{-1} \in C^{(a)}_{i,j} \}= \mathrm{Card} \{ (z_1, z_2) : I \eta^{(d)}_{p, j+l} = I z_1 z_2 \}. \] The coefficient $c(\alpha, \beta; \sigma)$ is finite for any $\sigma$ by Theorem \ref{thm-product}.  Similarly, we write \[ \sigma * \gamma = \sum_\tau c(\sigma, \gamma; \tau) \tau, \text{ where } \tau = \chi^{(e)}_{r, j+l+n}, \ e=1, 2, \ r \in \mathbb{Z}. \] We define \[ c(\alpha, \beta, \gamma; \tau)= \mathrm{Card} \{ (z_1, z_2, z_3) : I \eta^{(e)}_{r, j+l+n} = I z_1 z_2 z_3 \}. \] Since $(\alpha * \beta)*\gamma$ is defined, the number $\sum_\sigma c(\alpha, \beta;\sigma) c(\sigma, \gamma; \tau)$ is finite. Now it is not difficult to see that \[ \sum_\sigma c(\alpha, \beta;\sigma) c(\sigma, \gamma; \tau) = c(\alpha, \beta, \gamma; \tau). \]
Similarly, one can show that \[ \sum_{\sigma'} c(\alpha, \sigma';\tau) c(\beta, \gamma; \sigma') = c(\alpha, \beta, \gamma; \tau),  \] where $\sigma'$ is defined with regard to $\alpha *(\beta * \gamma)$. It proves the assertion of the proposition.
\end{proof}

\begin{Rmk}
The argument is essentially the same as in the case of Hecke operators on the space of modular forms; see \cite{Bu, Sh}.
\end{Rmk}

\vskip 1 cm

\section{The structure of $\mathcal{H}(G,I)$}

In this section, we investigate the structure of the Hecke algebra $\mathcal{H}(G,I)$. We will see that it has a big commutative subalgebra. After that, the center will be determined and Iwahori-Matsumoto type
relations will be found.

\vskip 5mm

For each $j \in \mathbb{Z}$, we let $\mathcal{H}_j$ be the subspace of $\mathcal{H}(G,I)$, consisting of the functions with their supports contained in $C_j$, where the set $C_j$ is defined in (\ref{eqn-cell}). We put \[ \mathcal{H}_-= \bigoplus_{j<0} \mathcal{H}_j \qquad \text{ and } \qquad \mathcal{H}_+= \bigoplus_{j>0} \mathcal{H}_j.\] Then, clearly, we have \[ \mathcal{H}(G,I)=\mathcal{H}_- \oplus \mathcal{H}_0 \oplus \mathcal{H}_+. \] It is easy to see that $\mathcal{H}_0$ is isomorphic to the (usual) affine Hecke algebra of $SL_2$ and that each $\mathcal{H}_j$ $(j \in \mathbb{Z})$ is a right and left $\mathcal{H}_0$-module.

We define \[\begin{array}{l} \Theta_{1,0}=\chi^{(1)}_{1,0} , \\
\Theta_{-1,0} = \chi^{(1)}_{-1,0} - (q-1) \chi^{(2)}_{-1,0} - (q-1) \chi^{(2)}_{0,0}  + q(q+q^{-1}-2) \chi^{(1)}_{0,0}, \\ \Theta_{0,1}=\chi^{(1)}_{0,1},\\ \Theta_{0,-1} = \chi^{(1)}_{0,-1} -(q-1) \displaystyle{\sum_{i\ge 0}} q^i \chi^{(2)}_{i,-1}.\end{array} \]
One can check $\Theta_{1,0}^{-1}= \Theta_{-1,0}$. (Recall that $\iota=q \chi^{(1)}_{0,0}$ is the identity.) The elements $\Theta_{1,0}$ and $\Theta_{-1,0}$ are the same as appear in Bernstein's presentation of the affine Hecke algebra $\mathcal{H}_0$.

\begin{Lem} \label{lem-relation}
The elements $\Theta_{1, 0}$, $\Theta_{-1, 0}$, $\Theta_{0, 1}$and $\Theta_{0,-1}$ commute with each other. \\ Moreover, we have:
\begin{enumerate} \item $\Theta_{1,0}^i*\Theta_{0,1}^j= q^{-(i+j-1)} \chi^{(1)}_{i,j}$ \quad for $i \in \mathbb{Z}$ and $j \ge 0$,

\item $\Theta_{-1,0}^i*\Theta_{0,-1}^j= q^{-(i+j-1)} \chi^{(1)}_{-i,-j} - (q-1) q^{-(i+j-1)} \chi^{(2)}_{-i,-j} + \displaystyle{\sum_{m>-i \atop a=1,2}} c_m^{(a)} \chi^{(a)}_{m, -j}$ \\ for $i \in \mathbb{Z}$ and $j > 0$ and for $c_m^{(a)} \in \mathbb{C}$.
\end{enumerate}
\end{Lem}

\begin{proof}
We use the formulas in Appendix and inductions on $i$ and $j$ to obtain (1) and (2).
\end{proof}

Let us denote by $\mathcal{A}$ the commutative subalgebra of $\mathcal{H}(G,I)$ generated by the elements $\Theta_{1, 0}$, $\Theta_{-1, 0}$, $\Theta_{0, 1}$and $\Theta_{0,-1}$. The structure of $\mathcal{A}$ is described by the following proposition.

\begin{Prop}
The algebra $\mathcal{A}$ is isomorphic to the quotient of the algebra $\mathbb{C}[X, X^{-1}, Y, Z]$ by the relation $YZ=0$. 
\end{Prop}

\begin{proof}
We have the surjective homomorphism $\phi : \mathbb{C}[X, X^{-1}, Y, Z]/(YZ) \rightarrow \mathcal{A}$ defined by \[ \phi (X) = \Theta_{1,0}, \quad \phi (X^{-1}) = \Theta_{-1,0}, \quad \phi (Y) = \Theta_{0,1} \quad \text{and} \quad \phi (Z) = \Theta_{0,-1}.\] It follows from (1) and (2) of Lemma \ref{lem-relation} that $\phi$ is also injective.
\end{proof}

\vskip 5mm

Our next task is to determine the center of $\mathcal{H}(G,I)$.

\begin{Thm}
The center of $\mathcal{H}(G,I)$ is the same as the center of $\mathcal{H}_0$ generated by the element $\Theta_{1,0}+\Theta_{-1,0}$.
\end{Thm}

\begin{proof}
It is well known that the element $\Theta_{1,0}+\Theta_{-1,0}$ generates the center of $\mathcal{H}_0$ (\cite{Lu}). Let us check if it commutes with $\chi^{(a)}_{i,j}$, $j \neq 0$. First, we assume $j>0$. Since $\chi^{(1)}_{i,j} \in \mathcal{A}$, it commutes with $\Theta_{1,0}+\Theta_{-1,0}$.  We have $\chi^{(2)}_{i,j}=q \chi^{(1)}_{i,j}*\chi^{(2)}_{0,0}$. Since $\chi^{(2)}_{0,0} \in \mathcal{H}_0$, we get $\chi^{(2)}_{i,j}*\left ( \Theta_{1,0}+\Theta_{-1,0} \right ) =  \left ( \Theta_{1,0}+\Theta_{-1,0} \right ) * \chi^{(2)}_{i,j}$. Next, we assume $j<0$. We obtain, using the formulas in Appendix, \[ \chi^{(1)}_{i,j}*\left ( \Theta_{1,0}+\Theta_{-1,0} \right ) =  q \chi^{(1)}_{i+1,j}+q^{-1} \chi^{(1)}_{i-1,j}= \left ( \Theta_{1,0}+\Theta_{-1,0} \right ) * \chi^{(1)}_{i,j}.\] Since $\chi^{(2)}_{i,j}=\chi^{(1)}_{i,j}*\chi^{(2)}_{0,0}-(1-q^{-1}) \chi^{(1)}_{i,j}$, the element  $\chi^{(2)}_{i,j}$ also commutes with $\Theta_{1,0}+\Theta_{-1,0}$. Therefore, the element $\Theta_{1,0}+\Theta_{-1,0}$ is in the center of the algebra $\mathcal{H}(G,I)$.

Suppose that $\zeta$ is an element in the center of $\mathcal{H}(G,I)$. We write $\zeta=\sum_{j \in \mathbb{Z}} \zeta_j$, $\zeta_j \in \mathcal{H}_j$. Since $\chi^{(2)}_{0,0} * \mathcal{H}_j \subset \mathcal{H}_j$ and $\mathcal{H}_j *\chi^{(2)}_{0,0} \subset \mathcal{H}_j$ for each $j$, the equality $\zeta * \chi^{(2)}_{0,0} = \chi^{(2)}_{0,0}*\zeta$ yields $\zeta_j * \chi^{(2)}_{0,0} = \chi^{(2)}_{0,0}*\zeta_j$ for each $j$. First, we assume $j>0$. Suppose that we choose the largest $i$ so that  \[\zeta_j = c_1 \chi^{(1)}_{i,j}+ c_2 \chi^{(2)}_{i,j} + \sum_{m<i \atop a=1,2} c_m^{(a)} \chi^{(a)}_{m,j}, \quad c_1 \neq 0 \text{ or } c_2 \neq 0.\] We get \[ \chi^{(2)}_{0,0}* \zeta_j = c_1(1-q^{-1}) \chi^{(1)}_{i,j} + c_2 (1-q^{-1}) \chi^{(2)}_{i,j} + \sum_{m<i \atop a=1,2} {c'}_m^{(a)} \chi^{(a)}_{m,j}.\] On the other hand, \[ \zeta_j*\chi^{(2)}_{0,0}= c_1 q^{-1} \chi^{(2)}_{i,j} + c_2 \chi^{(1)}_{i,j}+  c_2 (1-q^{-1}) \chi^{(2)}_{i,j} + \sum_{m<i \atop a=1,2} {c''}_m^{(a)} \chi^{(a)}_{m,j}.\] 
Thus we have $c_1=c_2=0$, a contradiction. It implies that $\zeta_j=0$. 

A similar argument also works for the case $j<0$, and we have  $\zeta_j=0$ in this case, too. Thus $\zeta =\zeta_0 \in \mathcal{H}_0$. It completes the proof.
\end{proof}

\vskip 5mm

The double affine Weyl group $W$ is not a Coxeter group, but it has a similar
presentation as one can see in the following proposition due to
A.\,N. Parshin. It is also related to Kyoji Saito's elliptic Weyl groups \cite{SaTa}.

\begin{Prop} \cite{Pa}
The group $W$ has a presentation
given by \begin{equation} \label{eqn-gr-relation}
 W \cong < s_0,s_1,s_2\ |\ s_0^2= s_1^2=
s_2^2=e,\
(s_0s_1s_2)^2 = e >.
\end{equation}
\end{Prop}

We can easily determine elements of $W$ corresponding to
the generators in the presentation. For example, we
can take, using the
same notation,
\[ s_0 = \begin{pmatrix} 0 & 1 \\ -1 & 0
\end{pmatrix}, \qquad s_1= \begin{pmatrix} 0 & t_1^{-1}
\\ -t_1 & 0
\end{pmatrix}, \qquad \text{and} \qquad
s_2= \begin{pmatrix} 0 & t_2^{-1} \\ -t_2 & 0
\end{pmatrix}.\]

We define  \begin{equation} 
\phi_0 =
q^{\frac 1 2}\chi^{(2)}_{0,0}= q^{\frac 1 2}\chi_{Is_0I} , \quad  \phi_1 =
q^{\frac 1 2}\chi^{(2)}_{-1,0}=q^{\frac 1 2} \chi_{Is_1I},
  \quad \text{and}
\quad  \phi_2 = q^{\frac 1 2}\chi^{(2)}_{0,-1}
=q^{\frac 1 2}\chi_{Is_2I}.
\end{equation}
The elements $\phi_0$ and $\phi_1$ have the special property
\[ \phi_0 * \mathcal{H}_-=0 \qquad \text{and} \qquad \phi_1*\mathcal{H}_+=0,\]
which follows from the formulas (2)(f) in Appendix. 


\begin{Prop}
The following identities
hold in $\mathcal{H}(G,I)$:
\begin{equation} \label{eqn-al-relation}
\begin{array}{ll}
& \phi_0*\phi_0 = (q^{\frac{1}{2}} - q^{-\frac{1}{2}})
\phi_0 + \iota, \\
& \phi_1*\phi_1 = (q^{\frac{1}{2}} - q^{-\frac{1}{2}})
\phi_1 + \iota, \\
\text{and }& \phi_0 *
\phi_1 * \phi_2 * \phi_0 * \phi_1 = \phi_2.
\end{array}
\end{equation}
\end{Prop}

\begin{proof}
We check all the relations using the formulas in Appendix.
\end{proof}

\begin{Rmk}
We can consider the relations in (\ref{eqn-al-relation}) as Iwahori-Matsumoto type relations. 
The first two relations in (\ref{eqn-al-relation}) are the usual deformation. The last one in (\ref{eqn-al-relation}) reflects the structure of the group algebra
of $W$. However, we have \[\phi_2*\phi_2 = (q^{\frac{1}{2}}
- q^{-\frac{1}{2}})
\sum_{m>0} q^{m-\frac 1 2} \chi^{(2)}_{m,-2},\] which reveals a new feature of the Hecke algebra $\mathcal{H}(G,I)$.
\end{Rmk}

\pagebreak

\section*{Appendix}

\begin{enumerate}
\item

\begin{enumerate}

\item If $(i,j) \ge (0,0)$ and $(k,l) \ge (0,0)$, or if
$(i,j) < (0,0)$ and $(k,l) < (0,0)$, then
\[ \chi^{(1)}_{i,j} * \chi^{(1)}_{k,l}
= q^{-1} \chi^{(1)}_{i+k,j+l} \quad \text{and} \quad
\chi^{(1)}_{i,j} * \chi^{(2)}_{k,l}
= q^{-1} \chi^{(2)}_{i+k,j+l}.
\]

\item If $i \ge 0$, $j=0$ and $l<0$, or if
$i < 0$, $j=0$ and $l>0$, then
\[ \chi^{(1)}_{i,0} * \chi^{(1)}_{k,l}
= q^{2 |i|-1} \chi^{(1)}_{i+k,l} \quad \text{and} \quad
\chi^{(1)}_{i,0} * \chi^{(2)}_{k,l}
= q^{2|i|-1} \chi^{(2)}_{i+k,l}.
\]

\item If $j>0$, $k<0$ and $l=0$, then
\[ \begin{array}{ll}  & \chi^{(1)}_{i,j} * \chi^{(1)}_{k,0}
= q^{-2k-1} \chi^{(1)}_{i+k,j}+(1-q^{-1})
\displaystyle{\sum_{m=i+k}^{i-k-1}}
 q^{i-k-m-1} \chi^{(2)}_{m,j}
 \\ \text{and} \quad &
\chi^{(1)}_{i,j} * \chi^{(2)}_{k,0}
=(1- q^{-1})
\displaystyle{\sum_{m=i+k+1}^{i-k-1}}
 q^{i-k-m-1} \chi^{(1)}_{m,j}+q^{-2k-2} \chi^{(2)}_{i+k,j}.
\end{array}
\]

\item If $j<0$, $k \ge 0$ and $l=0$, then
\[ \begin{array}{ll}  & \chi^{(1)}_{i,j} * \chi^{(1)}_{k,0}
= q^{2k-1} \chi^{(1)}_{i+k,j}+(1-q^{-1})
\displaystyle{\sum_{m=i-k}^{i+k-1}}
 q^{-i+k+m} \chi^{(2)}_{m,j}
 \\ \text{and} \quad &
\chi^{(1)}_{i,j} * \chi^{(2)}_{k,0}
=(1- q^{-1})
\displaystyle{\sum_{m=i-k}^{i+k}}
 q^{-i+k+m} \chi^{(1)}_{m,j}+q^{2k} \chi^{(2)}_{i+k,j}.
\end{array}
\]

\item If $i \ge 0$, $j=0$, $k<0$ and $l=0$, then
\[ \begin{array}{ll}  & \chi^{(1)}_{i,0} * \chi^{(1)}_{k,0}
= q^{\min\{2i-1,-2k-1 \}}
\chi^{(1)}_{i+k,0}+(1-q^{-1})
\displaystyle{\sum_{m=\max \{i+k,-i-k \}}^{i-k-1}}
 q^{i-k-m-1} \chi^{(2)}_{m,0}
 \\ \text{and} \quad &
\chi^{(1)}_{i,0} * \chi^{(2)}_{k,0}
=(1- q^{-1})
\displaystyle{\sum_{m=\max\{i+k+1,-i-k\}}^{i-k-1}}
 q^{i-k-m-1} \chi^{(1)}_{m,0}+q^{\min\{2i-1, -2k-2 \}}
 \chi^{(2)}_{i+k,0}.
\end{array}
\]

\item If $i<0$, $j=0$, $k \ge 0$ and $l=0$, then
\[ \begin{array}{ll}  & \chi^{(1)}_{i,0} * \chi^{(1)}_{k,0}
= q^{\min\{-2i-1,2k-1\}} \chi^{(1)}_{i+k,0}+(1-q^{-1})
\displaystyle{\sum_{m=i-k}^{\min\{i+k-1,-i-k-1\} }}
 q^{-i+k+m} \chi^{(2)}_{m,0}
 \\ \text{and} \quad &
\chi^{(1)}_{i,0} * \chi^{(2)}_{k,0}
=(1- q^{-1})
\displaystyle{\sum_{m=i-k}^{\min\{i+k,-i-k-1\} }}
 q^{-i+k+m} \chi^{(1)}_{m,0}+q^{\min\{-2i-1,2k\}}
 \chi^{(2)}_{i+k,0}.
\end{array}
\]

\end{enumerate}

\item

\begin{enumerate}

\item If $j>0$ and $l>0$, then
\[ \begin{array}{ll}  & \chi^{(2)}_{i,j} * \chi^{(1)}_{k,l}
= (1-q^{-1})
\displaystyle{\sum_{m \le i+k} }
 q^{i+k-m} \chi^{(1)}_{m,j+l}
 \\ \text{and} \quad &
\chi^{(2)}_{i,j} * \chi^{(2)}_{k,l}
= (1-q^{-1})
\displaystyle{\sum_{m \le i+k} }
 q^{i+k-m} \chi^{(2)}_{m,j+l}.
\end{array}
\]

\item If $j<0$ and $l<0$, then
\[ \begin{array}{ll}  & \chi^{(2)}_{i,j} * \chi^{(1)}_{k,l}
= (1-q^{-1})
\displaystyle{\sum_{m > i+k} }
 q^{-i-k+m-1} \chi^{(1)}_{m,j+l}
 \\ \text{and} \quad &
\chi^{(2)}_{i,j} * \chi^{(2)}_{k,l}
= (1-q^{-1})
\displaystyle{\sum_{m > i+k} }
 q^{-i-k+m-1} \chi^{(2)}_{m,j+l}.
\end{array}
\]

\item If $(i,j) \ge (0,0)$, $k<0$ and $l=0$,
or if $(i,j) < (0,0)$, $k \ge 0$ and $l=0$, then
\[ \chi^{(2)}_{i,j}*\chi^{(1)}_{k,0}
=q^{-1} \chi^{(2)}_{i-k,j} \quad \text{and} \quad
   \chi^{(2)}_{i,j}*\chi^{(2)}_{k,0}
=q^{-1} \chi^{(1)}_{i-k,j}.\]

\item IF $j>0$, $k \ge 0$ and $l=0$, then
\[ \begin{array}{ll}  & \chi^{(2)}_{i,j} * \chi^{(1)}_{k,0}
= (1-q^{-1})
\displaystyle{\sum_{m = i-k+1}^{i+k} }
 q^{i+k-m} \chi^{(1)}_{m,j} + q^{2k-1} \chi^{(2)}_{i-k,j}
 \\ \text{and} \quad &
\chi^{(2)}_{i,j} * \chi^{(2)}_{k,0}
= q^{2k} \chi^{(1)}_{i-k,j}+
(1-q^{-1})
\displaystyle{\sum_{m = i-k}^{i+k} }
 q^{i+k-m} \chi^{(2)}_{m,j}.
\end{array}
\]

\item If $j<0$, $k<0$ and $l=0$, then

\[ \begin{array}{ll}  & \chi^{(2)}_{i,j} * \chi^{(1)}_{k,0}
= (1-q^{-1})
\displaystyle{\sum_{m = i+k+1}^{i-k} }
 q^{-i-k+m-1} \chi^{(1)}_{m,j} + q^{-2k-1}
 \chi^{(2)}_{i-k,j}
 \\ \text{and} \quad &
\chi^{(2)}_{i,j} * \chi^{(2)}_{k,0}
= q^{-2k-2} \chi^{(1)}_{i-k,j}+
(1-q^{-1})
\displaystyle{\sum_{m = i+k+1}^{i-k-1} }
 q^{-i-k+m-1} \chi^{(2)}_{m,j}.
\end{array}
\]

\item If $i \ge 0$, $j=0$ and $l<0$, or if
$i < 0$, $j=0$ and $l>0$, then

\[ \chi^{(2)}_{i,0}*\chi^{(1)}_{k,l}=0 \quad
\text{and} \quad   \chi^{(2)}_{i,0}*\chi^{(2)}_{k,l}=0.\]

\item If $i \ge 0$, $j=0$ and $l>0$, then

\[ \begin{array}{ll}  & \chi^{(2)}_{i,0} * \chi^{(1)}_{k,l}
= (1-q^{-1})
\displaystyle{\sum_{m = -i+k}^{i+k} }
 q^{i+k-m} \chi^{(1)}_{m,l}
 \\ \text{and} \quad &
\chi^{(2)}_{i,0} * \chi^{(2)}_{k,l}
= (1-q^{-1})
\displaystyle{\sum_{m = -i+k}^{i+k} }
 q^{i+k-m} \chi^{(2)}_{m,l}.
\end{array}
\]

\item If $i<0$, $j=0$ and $l<0$, then

\[ \begin{array}{ll}  & \chi^{(2)}_{i,0} * \chi^{(1)}_{k,l}
= (1-q^{-1})
\displaystyle{\sum_{m = i+k+1}^{-i+k-1} }
 q^{-i-k+m-1} \chi^{(1)}_{m,l}
 \\ \text{and} \quad &
\chi^{(2)}_{i,0} * \chi^{(2)}_{k,l}
= (1-q^{-1})
\displaystyle{\sum_{m = i+k+1}^{-i+k-1} }
 q^{-i-k+m-1} \chi^{(2)}_{m,l}.
\end{array}
\]

\item If $i \ge 0$, $j=0$, $k \ge 0$ and $l=0$, then

\[ \begin{array}{ll}  & \chi^{(2)}_{i,0} * \chi^{(1)}_{k,0}
= (1-q^{-1})
\displaystyle{\sum_{m = \max\{i-k+1,-i+k\}}^{i+k} }
 q^{i+k-m} \chi^{(1)}_{m,0} + q^{\min\{2i,2k-1\}}
 \chi^{(2)}_{i-k,0}
 \\ \text{and} \quad &
\chi^{(2)}_{i,0} * \chi^{(2)}_{k,0}
= q^{\min\{2i,2k\}} \chi^{(1)}_{i-k,0}+
(1-q^{-1})
\displaystyle{\sum_{m =\max\{i-k, -i+k\}}^{i+k} }
 q^{i+k-m} \chi^{(2)}_{m,0}.
\end{array}
\]

\item If $i <0$, $j=0$, $k<0$ and $l=0$, then

\[ \begin{array}{ll}  & \chi^{(2)}_{i,0} * \chi^{(1)}_{k,0}
= (1-q^{-1})
\displaystyle{\sum_{m = i+k+1}^{\min\{i-k,-i+k-1\} } }
 q^{-i-k+m-1} \chi^{(1)}_{m,0} + q^{\min\{-2i-2,-2k-1\} }
 \chi^{(2)}_{i-k,0}
 \\ \text{and} \quad &
\chi^{(2)}_{i,0} * \chi^{(2)}_{k,0}
= q^{\min\{-2i-2,-2k-2\} } \chi^{(1)}_{i-k,0}+
(1-q^{-1})
\displaystyle{\sum_{m = i+k+1}^{\min\{-i+k-1,i-k-1\} } }
 q^{-i-k+m-1} \chi^{(2)}_{m,0}.
\end{array}
\]

\end{enumerate}

\end{enumerate}

\newpage


\end{document}